\documentclass[a4paper,10pt,leqno]{article}
\usepackage{amssymb,latexsym}
\usepackage{mltex}
\usepackage{amsmath}
 \usepackage{lscape}
\usepackage{enumerate}
\usepackage{amsthm}
\usepackage{hyperref}
\newtheorem{thm}{Theorem}

\newtheorem{lem}{Lemma}[section]

\newtheorem{rem}{Remark}

\newtheorem*{thrm}{Theorem}
\newcommand{\nablab}{\overline{\nabla}}

\newcommand{\iid}{\mathrm{Id}\,}

\newcommand{\ddet}{\mathrm{det}\,}

\newcommand{\trace}{\mathrm{tr\,}}

\newcommand{\C}{\mathbb{C}}

\newcommand{\Ss}{\mathbb{S}}
\newcommand{\HH}{\mathbb{H}}

\newcommand{\R}{\mathbb{R}}

\newcommand{\M}{\mathbb{M}}

\newcommand{\pre}{\Re e}

\newcommand{\beqt}{\begin{equation}}  \newcommand{\eeqt}{\end{equation}}
\newcommand{\bal}{\begin{align}}      \newcommand{\eal}{\end{align}}
\newcommand{\ba}{\begin{array}}      \newcommand{\ea}{\end{array}}
\newcommand{\bc}{\begin{center}}     \newcommand{\ec}{\end{center}}
\newcommand{\be}{\begin{enumerate}}  \newcommand{\ee}{\end{enumerate}}
\newcommand{\beq}{\begin{eqnarray}}  \newcommand{\eeq}{\end{eqnarray}}
\newcommand{\beQ}{\begin{eqnarray*}} \newcommand{\eeQ}{\end{eqnarray*}}
\newcommand{\bi}{\begin{itemize}}    \newcommand{\ei}{\end{itemize}}
\newcommand{\bt}{\begin{tabular}}    \newcommand{\et}{\end{tabular}}
\newcommand{\finpreuve}{\hfill\square\\}

\def\pf{\noindent{\textit {Proof :} }}

\title{Spinorial Characterizations of Surfaces into 3-dimensional pseudo-Riemannian Space Forms}
\author{Marie-Am\'elie Lawn\footnote{marie-amelie.lawn@uni.lu}\;\, and Julien Roth\footnote{julien.rothoth@univ-mlv.fr}
}
\date{}
\begin{document}
\maketitle

\begin{abstract}
We give a spinorial characterization of isometrically immersed
surfaces of arbitrary signature into 3-dimensional pseudo-Riemannian
space forms. For Lorentzian surfaces, this generalizes a recent work
of the first author in $\mathbb{R}^{2,1}$ to other Lorentzian space
forms. We also characterize immersions of Riemannian surfaces in
these spaces. From this we can deduce analogous results for timelike
immersions of Lorentzian surfaces in space forms of corresponding
signature, as well as for spacelike and timelike immersions of
surfaces of signature (0,2), hence achieving a complete spinorial
description for this class of pseudo-Riemannian immersions.
% We give a spinorial characterization of
%isometrically immersed surfaces of arbitrary signature into
%3-dimensional pseudo-Riemannian space forms. For Lorentzian
%surfaces, this generalizes a recent work of the first author in the
%3-dimensional Minkowski space to other Lorentzian space forms.
\end{abstract}
%%%%%%%%%%%%%%%%%%%%%%%%%%%%%%%%%%%%%%%%%%%%%%%%%%%%%%%%%%%%%%%%%%%%%%%%%%%%%%%%%%%
{\it keywords:} Dirac Operator, Killing Spinors, Isometric Immersions, Gauss and Codazzi Equations.\\\\
\noindent
{\it subjclass:} Differential Geometry, Global Analysis, 53C27, 53C40,
53C80, 58C40.

\date{}
%%%%%%%%%%%%%%%%%%%%%%%%%%%%%%%%%%%%%%%%%%%%%%%%%%%%%%%%%%%%%%%%%%%%%%%%%%%%%%%%%%%
\maketitle\pagenumbering{arabic}
%%%%%%%%%%%%%%%%%%%%%%%%%%%%%%%%%%%%%%%%%%%%%%%%%%%%%%%%%%%%%%%%%%%%%%%%%%%%%%%%%%%

\section{Introduction}
A fundamental question in the theory of submanifolds is to know
whether a (pseudo-)Riemannian manifold $(M^{p,q},g)$ can be
isometrically immersed into a fixed ambient manifold
$(\overline{M}^{r,s},\overline{g})$. In this paper, we focus on the
case of hypersurfaces, that is, codimension $1$. When the ambient
space is a space form, as the pseudo-Euclidean space $\R^{p,q}$ and
the pseudo-spheres $\Ss^{p,q}$ of positive constant curvature, or
the pseudo-hyperbolic spaces $\HH^{p,q}$ of negative constant
curvature, the answer is given by the well-known fundamental theorem
of hypersurfaces:
\begin{thrm}{\cite{On}}\label{thmOn} $(M^{p,q},g)$ be a pseudo-Riemannian manifold with signature $(p,q)$, $p+q=n$. Let $A$ be a symmetric Codazzi tensor , that is, $d^{\nabla}A=0$, satisying
\beQ
R(X,Y)Z= \delta\Big[ \left\langle A(Y),Z\right\rangle A(X)-\left\langle A(X),Z\right\rangle A(Y)\Big]
+\kappa\Big[ \left\langle Y,Z\right\rangle X-\left\langle X,Z\right\rangle Y\Big]
\eeQ
with $\kappa\in\R$ for all $x\in M$ and $X,Y,Z\in T_xM$.\\
Then, if $\delta=1$ (resp. $\delta=-1$), there exists locally a
spacelike (resp. timelike) isometric immersion of $M$ in
$\M^{p+1,q}(\kappa)$ (resp. $\M^{p,q+1}(\kappa)$).
\end{thrm}
In the Riemannian case and for small dimensions ($n=2$ or $3$), an other necessary and
sufficient condition is now well-known. This condition is expressed
in spinorial terms, namely, by the existence of a special spinor field.
This work initiated by Friedrich \cite{Fr} in the late 90's for surfaces of $\R^3$
was generalized for surfaces of $\Ss^3$ and $\HH^3$ \cite{Mo} and other 3-dimensional homogeneous manifolds \cite{Roth4}.\\

\indent The first author \cite{La} uses this approach to give a
spinorial characterization of space-like immersions of Lorentzian
surfaces in the Minkowski space $\R^{2,1}$. In this paper, we give a
generalization of this result to Lorentzian or Riemannian surfaces
into one of the three Lorentzian space forms, $\R^{2,1}$, $\Ss^{2,1}$
or $\widetilde{\HH}^{2,1}$. These finally allows us to give a
complete spinorial characterization for spacelike as well as for
timelike immersions of surfaces of arbitrary signature into
pseudo-Riemannian space forms.

\indent We will begin by a section of recalls about extrinsic
pseudo-Riemannian spin geometry. For further details, one refers to
\cite{BGM,Bau} for  basic facts about spin geometry and
\cite{BGM,BM,LM} for the extrinsic aspect.

\section{Preliminaries}
\subsection{Pseudo-Riemannian spin geometry}
\label{sec21} Let $(M^{p,q},g)$, $p+q=2$, be an oriented
pseudo-Riemannian surface of arbitrary signature isometrically
immersed into a three-dimensional pseudo-Riemannian spin manifold
$(N^{r,s},\overline{g})$. We introduce the parameter $\varepsilon$ as
follows: $\varepsilon=i$ if the immersion is timelike and
$\varepsilon=1$ if the immersion is spacelike. Let $\nu$ be a unit
vector normal to $M$. The fact that $M$ is oriented implies that $M$
carries a spin structure induced from the spin structure of $N$ and
we have the following identification of the spinor bundles and
Clifford multiplications:
$$\left\lbrace\begin{array}{ll}
\Sigma N_{|M}\equiv\Sigma M.\\
X\cdot \varphi_{|M}=\big(\varepsilon\nu\bullet X\bullet\varphi)_{|M},
\end{array}\right.$$
where $\cdot$ and $\bullet$ are the Clifford multiplications, respectively on $M$ and $N$. Moreover, we have the following well-known spinorial Gauss formula
\beqt\label{gauss}
\nablab_X\varphi=\nabla_X\varphi-\frac{\varepsilon}{2}A(X)\cdot\varphi,
\eeqt
where $A$ is the shape operator of the immersion. Finally, we recall the Ricci identity on $M$
\beqt\label{ricci}
R(e_1,e_2)\varphi=\frac{1}{2}\varepsilon_1\varepsilon_2R_{1221}e_1\cdot e_2\cdot\varphi,
\eeqt
 where ${e_1,e_2}$ is a local orthonormal frame of $M$ and $\varepsilon_j=g(e_j,e_j)$.\\
 \noindent
 The complex volume element on the surface depends on the signature and is defined by
$$\omega_{p,q}^{\mathbb{C}}=
    i^{q+1}e_1\cdot e_2.$$
Obviously ${\omega_{p,q}^{\mathbb{C}}}^2=1$ independently of the
signature and the action of $\omega^{\C}$ splits $\Sigma M$ into two
eigenspaces $\Sigma^{\pm}M$ of real dimension $2$. Therefore, a
spinor field $\varphi$ can be written as
$\varphi=\varphi^++\varphi^-$ with
$\omega^{\C}\cdot\varphi^{\pm}=\pm\varphi^{\pm}$. Finally, we denote
$\overline{\varphi}=\omega^{\C}\cdot\varphi=\varphi^+-\varphi^-$.
\subsection{Restricted Killing spinors}
Let $(M^{p,q},g)$, $p+q=2$ be a surface of the pseudo-Riemannian space form
$\M^{r,s}(\kappa)$, $r+s=3$, $p\leqslant r$, $q\leqslant s$. This space form carries a Killing spinor
$\varphi$, that is satisfying $\nablab_X\varphi=\lambda
X\bullet\varphi$, with $\kappa=4\lambda^2$. From the Gauss formula
$\eqref{gauss}$, the restriction of $\varphi$ on $M$ satisfies the equation
\beqt \label{Special_Killing}
\nabla_X\varphi=\frac{\varepsilon}{2}A(X)\cdot\varphi+\lambda
X\bullet\varphi.
\eeqt
But we have
$$X\bullet\varphi=\varepsilon^2\nu\bullet\nu\bullet
X\bullet\varphi=-\varepsilon^2\nu\bullet
X\bullet\nu\bullet\varphi=-\varepsilon X\cdot(\nu\bullet\varphi).$$
Moreover, the complex volume element
$\omega^{\mathbb{C}}_{r,s}=-i^se_1\bullet e_2\bullet \nu$ of
$\M^{r,s}(\kappa)$ over $M$ acts as the identity on $\Sigma
\M^{r,s}(\kappa)_{|M}\equiv \Sigma M$. Thus, we have \beQ
\nu\bullet\varphi&=&\omega^{\C}_{r,s}\bullet\varphi=-i^s\nu\bullet e_1\bullet e_2\bullet\nu\bullet\varphi\\
&=& i^s\nu \bullet e_1\bullet \nu \bullet e_2\bullet\varphi\\
&=& i^s\varepsilon^2 (\varepsilon\nu\bullet e_1)\bullet(\varepsilon \nu e_2\bullet\varphi)\\
&=&i^s\varepsilon^2e_1\cdot e_2\cdot\varphi.
\eeQ
Hence a simple case
by case computation shows that we have $$X\bullet\varphi
=-i^s\varepsilon^3X\cdot e_1\cdot
e_2\cdot\varphi=iX\cdot\omega^{\C}_{p,q}\cdot\varphi=iX\cdot\overline{\varphi}.$$
in the six possible cases (for the three possible signatures (2,0),
(1,1), (0,2) of the surface with respectively $\varepsilon=1$ or
$i$).\\
We will call a spinor solution of equation \eqref{Special_Killing} a real special Killing spinor (RSK)-spinor if $\varepsilon\in\mathbb{R}$,
and an imaginary special Killing spinor (ISK)-spinor if $\varepsilon\in i\mathbb{R}$.
\subsection{Norm assumptions}
In this section, we precise the norm assumptions useful for the statement of the main result. Let $(M^{p,q},g)$ be a pseudo-Riemannian surface and $\varphi$ spinor field on $M$. Let $\varepsilon=1$ or $i$ and $\lambda\in\R$ or $i\R$. We say that $\varphi$ satisfies the norm assumption $\mathcal{N}_{\pm}(p,q,\lambda,\varepsilon)$ if the following holds:
\begin{enumerate}
\item For $p=2,q=0$ or $p=0,q=2$:
\begin{itemize}
\item If $\varepsilon=1$, then $X|\varphi_1|^2=\pm2\pre\left\langle i\eta X\cdot\overline{\varphi},\varphi\right\rangle .$
\item If $\varepsilon=i$, then $X\left\langle\varphi,\overline{\varphi}\right\rangle=\pm2\pre\left\langle i\eta X\cdot\varphi_1,\varphi_1\right\rangle.$
\end{itemize}
\item For $p=1,q=1$: $\varphi$ is non-isotropic

\end{enumerate}
\section{The main result}
We now state the main result of the present paper.
\begin{thm}\label{thm1}
Let $(M^{p,q},g)$, $p+q=2$ be an oriented pseudo-Riemannian
manifold. Let $H$ be a real-valued function. Then, the three
following statements are equivalent:
\begin{itemize}
\item[1.] There exist two nowhere vanishing spinor fields $\varphi_1$ and $\varphi_2$ satisfying the norm assumptions
$\mathcal{N}_{+}(p,q,\lambda,\varepsilon)$ and $\mathcal{N}_{-}(p,q,\lambda,\varepsilon)$ respectively and
$$D\varphi_1=2\varepsilon
H\varphi_1+2i\lambda \overline{\varphi_1}\quad\text{and}\quad
D\varphi_2=-2\varepsilon H\varphi_2-2i\lambda\overline{\varphi}_2.$$
\item[2.] There exist two spinor fields $\varphi_1$ and $\varphi_2$ satisfying
$$\nabla_X\varphi_1=\frac{\varepsilon}{2} A(X)\cdot\varphi_1-i\lambda X\cdot\overline{\varphi}_1,\quad\text{and}
\quad\nabla_X\varphi_2=-\frac{\varepsilon}{2} A(X)\cdot\varphi_1+i\lambda X\cdot\overline{\varphi}_2,$$
where $A$ is a $g$-symmetric endomorphism and $H=-\frac{1}{2}\trace(A)$.

\item[3.] There exists a local isometric immersion from $M$ into the (pseudo)-Riemannian
space form $\M^{p+1,q}(4\lambda^2)$ ({\it resp.} $\M^{p,q+1}(4\lambda^2)$) if $\varepsilon=1$ ({\it resp.}  $\varepsilon=i$) with mean curvature $H$
and shape operator $A$.
\end{itemize}
\end{thm}
\begin{rem}
Note that, in this result, two spinor fields are needed to get an
isometric immersion. Nevertheless, for the case of Riemannian
surfaces in Riemannian space forms (Friedrich \cite{Fr} and Morel
\cite{Mo}) only one spinor solution of one of the two equations is
sufficient. This is also the case for surfaces of signature $(0,2)$
in space forms of signature $(0,3)$.
\end{rem}
\noindent
In order to prove this theorem, we give two technical lemmas.
\begin{lem}\label{lem1}
Let $(M^{p,q},g)$ be an oriented (pseudo)-Riemannian surface and
$\eta,\lambda$ two complex numbers. If $M$ carries a spinor field
satisfying
$$\nabla_X\varphi=\eta A(X)\cdot\varphi+i\lambda X\cdot\overline{\varphi},$$
then, we have
$$\left( \varepsilon_1\varepsilon_2R_{1212}+4\eta^2\ddet(A)-\lambda^2\right)e_1\cdot e_2\cdot\varphi=2\eta d^{\nabla}A(e_1,e_2)\cdot\varphi.$$
\end{lem}
\noindent
\pf An easy computation yields
\begin{eqnarray*}
\nabla_X\nabla_Y\varphi&=&\eta \nabla_XA(Y)\cdot\varphi+\eta^2A(Y)\cdot A(X)\cdot\varphi+i\eta\lambda  A(Y)\cdot X\cdot\omega^{\mathbb{C}}\cdot\varphi\\ \\
&&+i\lambda
\nabla_XY\cdot\omega^{\mathbb{C}}\cdot\varphi+i\eta\lambda
Y\cdot\omega^{\mathbb{C}}\cdot
A(X)\cdot\varphi-\lambda^2Y\cdot\omega^{\mathbb{C}}\cdot
X\cdot\omega^{\mathbb{C}}\varphi.
\end{eqnarray*}
Hence (the other terms vanish by symmetry)
\begin{eqnarray*}
\mathcal{R}(e_1,e_2)\varphi&=&\nabla_{e_1}\nabla_{e_2}\varphi-\nabla_{e_2}\nabla_{e_1}\varphi-\nabla_{[e_1,e_2]}\varphi\\
&=&\eta\big(\nabla_{e_1}A(e_2)-\nabla_{e_2}A(e_1)-A([e_1,e_2])\big)\varphi+\eta^2(A(e_2)A(e_1)-A(e_1)A(e_2))\varphi\\
&&-\lambda^2(e_2\cdot\omega^{\mathbb{C}}\cdot
e_1\omega^{\mathbb{C}}-e_1\cdot\omega^{\mathbb{C}}\cdot
e_2\omega^{\mathbb{C}})\varphi.
\end{eqnarray*}
Since we have
$$ A(e_2)A(e_1)-A(e_1)A(e_2)=-2\det(A)$$
 and
$$e_2\cdot\omega^{\mathbb{C}}\cdot e_1\omega^{\mathbb{C}}-e_1\cdot\omega^{\mathbb{C}}\cdot e_2\cdot\omega^{\mathbb{C}}=
e_1\cdot e_2\cdot\left(\omega^{\mathbb{C}}\right)^2-e_2\cdot e_1\cdot\left(\omega^{\mathbb{C}}\right)^2=2e_1\cdot e_2,$$
by the Ricci identity $\eqref{ricci}$, we get
\begin{eqnarray*}
\frac{1}{2}\varepsilon_1\varepsilon_2R_{1221}e_1e_2\cdot\varphi=\eta d^{\nabla}A(e_1,e_2)-2\eta^2\det(A)e_1\cdot e_2\varphi-2\lambda^2 e_1\cdot e_2 \cdot \varphi,
\end{eqnarray*}
and finally
\begin{eqnarray}\label{Gauss_codazzi}
(-\varepsilon_1\varepsilon_2 R_{1212}+4\eta^2\det(A)+4\lambda^2)e_1\cdot e_2\cdot\varphi=
2\eta d^{\nabla}A(e_1,e_2)\cdot\varphi.
\end{eqnarray}
$\finpreuve$
Now, we give this second lemma
\begin{lem}\label{lem2}
Let $(M^{p,q},g)$ be an oriented (pseudo)-Riemannian surface and
$\lambda$ a complex number. If $M$ carries a spinor field solution
of the  equation \begin{eqnarray}\label{Dirac_equation}D\varphi=\pm\left(\varepsilon
H\varphi+2i\lambda\overline{\varphi}\right)\end{eqnarray} and satisfying the norm
assumption $\mathcal{N}_{\pm}(p,q,\lambda,\varepsilon)$, then this spinor
satisfies
$$\nabla_X\varphi=\pm\left(\frac{\varepsilon}{2} A(X)\cdot\varphi-i\lambda X\cdot\overline{\varphi}\right).$$
\end{lem}
\noindent \pf We give the proof for the sign $+$. The other case is strictly the same .\\
{\bf Case of signature (1,1)}: We define the
endomorphism $B_{\varphi}$ by
$$(B_{\varphi})^i_j=g(B_{\varphi}(e_i),e_j)=\beta_{\varphi}(e_i,~e_j):=\langle\varepsilon\nabla_{e_i}\varphi,~e_j\cdot\varphi\rangle.$$
Obviously
%$tr(B_{\varphi})=g^{ij}{B_{ij}}_{\varphi}=\frac{1}{2}H-2i\lambda\omega^{\C}$.
 Using $\frac{e_i\cdot\varphi^{\pm}}{\langle \varphi^+,\varphi^-\rangle}$
as a normalized dual frame of $\Sigma^{\mp}M$ and the same proof as
in \cite{La} we can show that
\begin{eqnarray*}\label{derniere_minute} \langle\nabla_X\varphi,e_i\cdot\varphi^{\pm}\rangle=\langle\varepsilon\nabla_X\varphi,\varepsilon e_i\cdot\varphi^{\pm}\rangle=-\frac{1}{2\varepsilon\langle\varphi^+,\varphi^-\rangle}\langle
B_{\varphi}(X)\cdot\varphi,e_i\cdot\varphi^{\mp}\rangle.\end{eqnarray*}
and hence
$\nabla_X\varphi=-\frac{1}{2\varepsilon\langle\varphi^+,\varphi^-\rangle}
B_{\varphi}(X)\cdot\varphi$. Moreover
\begin{eqnarray*}
\beta_{\varphi}(e_1,~e_2)&=&\langle\nabla_{e_1}\varphi,~e_2\cdot\varphi\rangle
=-\langle\varepsilon\nabla_{e_1}\varphi,~e_1^2\cdot e_2\cdot\varphi\rangle\\
                  &=&-\langle\varepsilon e_1\cdot\nabla_{e_1}\varphi,~e_1\cdot e_2\cdot\varphi\rangle
                  =-\langle\varepsilon D\varphi+\varepsilon e_2\cdot\nabla_{e_2}\varphi,~e_1\cdot e_2\cdot\varphi\rangle\\
                  &=&-\varepsilon^2 H\langle\varphi,~e_1\cdot
                  e_2\cdot\varphi\rangle-\langle2i\varepsilon\lambda\omega^{\mathbb{C}}\cdot\varphi,~e_1\cdot
                  e_2\cdot\varphi\rangle +\beta_{\varphi}(e_2,e_1)\\
                  &=&-\langle2i\varepsilon\lambda\omega^{\mathbb{C}}\cdot\varphi,~e_1\cdot
                  e_2\cdot\varphi\rangle +\beta_{\varphi}(e_2,e_1),
\end{eqnarray*}
since for any $\varphi,\psi\in\Gamma(\Sigma M)$
  $$\langle\varphi,~e_1\cdot
e_2\cdot\psi\rangle=\langle e_2\cdot
e_1\cdot\varphi,~\psi\rangle=-\langle e_1\cdot
e_2\cdot\varphi,~\psi\rangle =-\langle \varphi,~e_1\cdot
e_2\cdot\psi\rangle=0.$$
Let now consider the decomposition $\beta_{\varphi}(X,Y)=S_{\varphi}(X,Y)+T_{\varphi}(X,Y)$
in the symmetric part $S_{\varphi}$ and antisymmetric part
$T_{\varphi}$. Hence, we see easily that if $\lambda/\varepsilon\in i\mathbb{R}$, then $\beta_{\varphi}$ is symmetric, {\it i.e.}, $T_{\varphi}=0$. and if $\lambda/\varepsilon\in\mathbb{R}$, then
$T_{\varphi}(X)=2i\lambda/\varepsilon\,\omega^{\mathbb{C}}\cdot X$. In the two cases, we have
$$ \nabla_X\varphi=\frac{\varepsilon}{2} A(X)\cdot\varphi-i\lambda X\cdot\overline{\varphi},$$
by setting $A=2S_{\varphi}$. We verify easily that $tr(A)=2tr(S_{\varphi})=2tr(B_{\varphi})=-2H$.\\ \\
{\bf Case of signature $(2,0)$ or $(0,2)$:} The proof is fairly
standard following the technique used in \cite{Fr}, \cite{Mo} and
\cite{Roth4}. We consider the tensors $Q^{\pm}_{\varphi}$ defined by
$$Q^{\pm}_{\varphi}(X,Y)=\pre\left\langle \varepsilon \nabla_X\varphi^{\pm},Y\cdot\varphi^{\mp}\right\rangle.$$
Then, we have
\beQ
\trace(Q^{\pm}_{\varphi})=-\pre\left\langle \varepsilon D\varphi^{\pm},\varphi^{\mp}\right\rangle=-\pre\left\langle \varepsilon(\varepsilon H\pm 2i\lambda\varphi^{\mp},\varphi^{\mp}\right\rangle=-\varepsilon^2\big(H\pm2\pre(\lambda)\big)|\varphi^{\mp}|^2.
\eeQ
Moreover, we have the following defect of symmetry of $Q^{\pm}_{\varphi}$,
\beQ
Q^{\pm}_{\varphi}(e_1,e_2)&=&\pre\left\langle \varepsilon\nabla_{e_1}\varphi^{\pm},e_2\cdot\varphi^{\mp}\right\rangle=\pre\left\langle \varepsilon e_1\cdot\nabla_{e_1}\varphi^{\pm},e_1\cdot e_2\cdot\varphi^{\mp}\right\rangle\\
&=&\pre\left\langle \varepsilon D\varphi^{\pm},e_1\cdot e_2\cdot\varphi^{\mp}\right\rangle-\pre\left\langle \varepsilon \nabla_{e_2}\varphi^{\pm},e_1\cdot e_2\cdot\varphi^{\mp}\right\rangle\\
&=&\pre\left\langle (\varepsilon^2H\pm 2i\varepsilon\lambda)\varphi^{\mp},e_1\cdot e_2\cdot\varphi^{\mp}\right\rangle+\pre\left\langle \varepsilon\nabla_{e_2}\varphi^{\pm},e_1\cdot\varphi^{\mp}\right\rangle\\
&=&2\pre(\varepsilon\lambda)|\varphi^{\mp}|^2+Q^{\pm}_{\varphi}(e_2,e_1).
\eeQ
Then, using the fact that
$\varepsilon e_1\cdot\dfrac{\varphi^{\pm}}{|\varphi^{\pm}|^2}$ and
$\varepsilon e_2\cdot\dfrac{\varphi^{\pm}}{|\varphi^{\pm}|^2}$ form a local
orthonormal frame of $\Sigma^{\mp}M$ for the real scalar product
$\pre\left\langle\cdot,\cdot\right\rangle$, we see easily that
$$\nabla_X\varphi^+=\varepsilon \frac{Q^+_{\varphi}(X)}{|\varphi^-|^2}\cdot\varphi^-\quad\text{and}\quad\nabla_X\varphi^-=\varepsilon \frac{Q^-_{\varphi}(X)}{|\varphi|^+|^2}\cdot\varphi^+.$$
We set $W=\frac{Q_{\varphi}^+}{|\varphi^-|^2}-\frac{Q_{\varphi}^-}{|\varphi^+|^2}$.
From the above computations, we have immediately that
$W+\pre\left(i\lambda/\varepsilon\right)\iid$ is symmetric and trace-free. Now, we will show that
$W+\pre\left(i\lambda/\varepsilon\right)\iid$ is of rank at most $1$. First, we have
$$X|\varphi^+|^2+\varepsilon^2X|\varphi^-|^2=2\pre\left\langle \varepsilon W(X)\cdot\varphi^-,\varphi^+\right\rangle.$$
Moreover, from the norm assuption $\mathcal{N}(p,q,\lambda,\varepsilon)$, we have
$$X|\varphi^+|^2+\varepsilon^2X|\varphi^-|^2=2
\pre\left\langle i \lambda X\cdot\varphi,\varphi\right\rangle=4
\pre\left\langle i \lambda X\cdot\varphi^-,\varphi^+\right\rangle.$$
 We
deduce immediately that $W+2\pre\left(i\lambda/\varepsilon\right)\iid$ is of rank at most $1$ and
hence vanishes identically since it is symmetric and trace-free.
Thus, we have the following relation
$$|\varphi^+|^2Q_{\varphi}^+-|\varphi^-|^2Q_{\varphi}^-=-2\pre(i\lambda/\varepsilon)|\varphi^+|^2|\varphi^-|^2g.$$
From now on, we will distinguish two cases.\\
$\bullet$ {\it Case 1:} $i\lambda/\varepsilon\in\R$.\\
Then we are in one of these two possible situations: $\varepsilon=i$ and $\lambda\in \R$ or $\varepsilon=1$ and $\lambda\in i\R$. The second situation was studied by Morel \cite{Mo}.\\
So we define the following tensor
$F:=Q_{\varphi}^+-Q_{\varphi}^-+2i\varepsilon\lambda(|\varphi^+|^2-|\varphi^-|^2)g$.
We have then \beQ
\nabla_X\varphi&=&\nabla_X\varphi^++\nabla_X\varphi^-=\varepsilon \frac{Q_{\varphi}^+(X)}{|\varphi^-|^2}\cdot\varphi^++\varepsilon \frac{Q_{\varphi}^-(X)}{|\varphi^-|^2}\cdot\varphi^-\\
&=& \varepsilon \frac{F(X)}{|\varphi|^2}\cdot(\varphi^++\varphi^-)-i\lambda X\cdot\varphi^--i\lambda X\cdot\varphi^+\\
&=&\frac{\varepsilon}{2}A(X)\cdot\varphi-i\lambda X\cdot\overline{\varphi},
\eeQ
where we have set $A=\frac{2F}{|\varphi^2|}$. We conclude by noticing that $A$ is a symmetric tensor with $\trace(A)=-2H$.\\
$\bullet$ {\it Case 2:} $i\lambda/\varepsilon\in i\R$.\\
Then we are in one of these two possible situations: $\varepsilon=i$ and $\lambda\in i\R$ or $\varepsilon=1$ and $\lambda\in \R$. The second situation was studied by Morel \cite{Mo}.\\
In this case, we have from the previous computations that $W$ vanishes identically. So we set
$$F=\frac{Q_{\varphi}^+}{|\varphi¬-|^2}=\frac{Q_{\varphi}^-}{|\varphi^+|^2}$$
 and then we have $\nabla_X\varphi=F(X)\cdot\varphi$, where $F(X)$ is defined by $g(F(X),Y)=F(X,Y)$. Nevertheless, $F$ is not symmetric. We define the following symmetric tensor $A(X,Y)=\frac{1}{|\varphi|^2}(F(X,Y)+F(Y,X))$. We compute immediately
$$A(e_1,e_1)=2F(e_1,e_1)/|\varphi|^2\quad,\quad A(e_2,e_2)=2F(e_2,e_2)/|\varphi|^2,$$
$$A(e_1,e_2)=2F(e_1,e_1)/|\varphi|^2-2\lambda/\varepsilon\quad\text{and}\quad A(e_2,e_2)=2F(e_2,e_2)/|\varphi|^2+2\lambda/\varepsilon.$$
Finally, we conclude that
$$\nabla_X\varphi=\frac{\varepsilon}{2}A(X)\cdot\varphi+\lambda X\cdot\omega\cdot\varphi=\frac{\varepsilon}{2}A(X)\cdot\varphi-i\lambda X\cdot\overline{\varphi}.$$
$\finpreuve$\\
Now, we can give the proof of Theorem \ref{thm1}.
We have already proven that 3. implies 2. which implies 1. Moreover, Lemma \ref{lem2} shows that 1. implies 2.
Now, we will prove that 2. implies 3. For this, we use Lemma \ref{lem1},
but we need to distinguish the three cases for the different signatures. Let $\varphi=\varphi^++\varphi^-$.\\\\
\noindent
\textbf{Case of signature (2,0)}: Here, $\omega^{\mathbb{C}}=ie_1e_2$, hence $e_1\cdot e_2\cdot\varphi=-i\omega^{\mathbb{C}}\cdot\varphi=-i\bar{\varphi}$.\\
    Hence formula \eqref{Gauss_codazzi} becomes
    \begin{eqnarray*}
-i\underbrace{(-R_{1212}+\varepsilon^2\det(A)+4\lambda^2)}_{G_{2,0}}\bar{\varphi}=\varepsilon\underbrace{d^{\nabla}A(e_1,e_2)}_{C_{2,0}}\cdot\varphi.
\end{eqnarray*}
or equivalently $\varepsilon C_{2,0}\cdot\varphi^{\pm}=\pm i G_{2,0}\varphi^{\mp}$. Applying two times this relation we have finally
$$ \varepsilon^2 \lvert\lvert C_{2,0}\lvert\lvert^2\varphi^{\pm}=-G^2_{2,0}\varphi^{\pm}.$$
Again we have two cases.\\ \noindent
$\bullet$ \textit{Spacelike immersion:} $\varepsilon=1$, $M^{2,0}\hookrightarrow\mathbb{M}^{3,0}$.\\
We refer to \cite{Fr} for the immersion in $\mathbb{R}^{3,0}$ and to
\cite{Mo} for $\mathbb{S}^{3}$ and $\mathbb{H}^{3}$. Only one (RSK)-spinor is needed.\\ \noindent
$\bullet$ \textit{Timelike immersion}: $\varepsilon=i$,
$M^{2,0}\hookrightarrow\mathbb{M}^{2,1}$.\\
Two (ISK)-spinors are needed. We deduce from the above relations between $\varphi^{\pm}_1$ and $\varphi^{\pm}_2$
that $\left\langle C_{2,0}\cdot\varphi_1,\varphi_2\right\rangle=0$. Moreover, in this case we have $\left\langle\varphi_1,\varphi_2\right\rangle=0$.
Thus, since the spinor bundle $\Sigma M$ is of complex rank $2$, we have $C_{2,0}\cdot\varphi_1=f\varphi_1$ where $f$ is a complex-valued function
over $M$. By taking the inner product by $\varphi_1$, we see immediately that $f$ only takes imaginary values, that is $f=ih$ with $h$ real-valued.
Thus, we have $\pm G_{2,0}\varphi^{\pm}_1=ih\varphi_1^{\pm}$. Since $\varphi_1^{\pm}$ do not vanish simultaneously, we deduce that $h$ and $G_{2,0}$
vanish identically. Thus $C$ vanishes too and the Gauss and Codazzi equation are satisfied. Then, we get the conclusion by the fundamental theorem
of hypersurfaces given above.\\\\ \noindent
 \textbf{Case of signature (1,1)}: $\omega^{\mathbb{C}}=-e_1e_2$, hence $e_1\cdot e_2\cdot\varphi=-\omega^{\mathbb{C}}\cdot\varphi=-\bar{\varphi}$.\\
    Hence formula \eqref{Gauss_codazzi} becomes
    \begin{eqnarray*}
-\underbrace{(R_{1212}+\varepsilon^2\det(A)+4\lambda^2)}_{G_{1,1}}\bar{\varphi}=\varepsilon\underbrace{d^{\nabla}A(e_1,e_2)}_{C_{1,1}}\cdot\varphi.
\end{eqnarray*}
or equivalently $\varepsilon C_{1,1}\cdot\varphi^{\pm}=G_{1,1}\varphi^{\mp}$. Applying two times this relation we have finally
$$ \varepsilon^2 \lvert\lvert C_{1,1}\lvert\lvert^2\varphi^{\pm}=G^2_{1,1}\varphi^{\pm}.$$
$\bullet$ \textit{Spacelike immersion:} $\varepsilon=1$, $M^{1,1}\hookrightarrow\mathbb{M}^{2,1}$.\\
We refer to \cite{La} for the immersion in $\mathbb{R}^{2,1}$. Let us consider the other space forms.\\
Here again, we need two (RSK)-spinors. Since $\varphi_1^{\pm}$ do not vanish at the same point, we have clearly that $||C_{1,1}||=G_{1,1}^2\geqslant0$. Moreover, we have
\beQ-||C_{1,1}||^2\left\langle\varphi_1,\varphi_2\right\rangle&=&\left\langle C_{1,1}\cdot\varphi_1,C_{1,1}\cdot\varphi_2\right\rangle\\
&=&-G_{1,1}^2\left\langle e_1\cdot e_2\varphi_1,e_1\cdot e_2\cdot\varphi_2\right\rangle\\
&=&G_{1,1}^2\left\langle\varphi_1,\varphi_2\right\rangle.
\eeQ
Since $\left\langle\varphi_1,\varphi_2\right\rangle$ never vanishes, we deduce that $||C_{1,1}||=-G_{1,1}^2\leqslant0$. Consequently, $||C_{1,1}||=G_{1,1}=0$. Moreover, $C_{1,1}$ is not isotropic. Indeed, since $G_{1,1}=0$, we have $C_{1,1}\cdot\varphi_1=0$ and thus $C_{1,1}$ automatically vanishes as proved in \cite{La}.\\
\noindent
$\bullet$ \textit{Timelike immersion}: $\varepsilon=i$,
$M^{1,1}\hookrightarrow\mathbb{M}^{1,2}$. It is easy to see that computations similar to the one for the previous case give the result.\\\
 Two (ISK)-spinors are needed.\\\\
 \noindent
\textbf{Case of signature (0,2)} $\omega^{\mathbb{C}}=-ie_1e_2$, hence $e_1\cdot e_2\cdot\varphi=i\omega^{\mathbb{C}}\cdot\varphi=i\bar{\varphi}$.\\
    Hence formula \eqref{Gauss_codazzi} becomes
    \begin{eqnarray*}
i\underbrace{(-R_{1212}+\varepsilon^2\det(A)+4\lambda^2)}_{G_{0,2}}\bar{\varphi}=\varepsilon\underbrace{d^{\nabla}A(e_1,e_2)}_{C_{0,2}}\cdot\varphi.
\end{eqnarray*}
or equivalently $\varepsilon C_{0,2}\cdot\varphi^{\pm}=\pm i G_{0,2}\varphi^{\mp}$. Applying two times this relation we have finally
$$ \varepsilon^2 \lvert\lvert C_{0,2}\lvert\lvert^2\varphi^{\pm}=-G^2_{0,2}\varphi^{\pm}.$$
$\bullet$ \textit{Spacelike immersion:} $\varepsilon=1$, $M^{0,2}\hookrightarrow\mathbb{M}^{1,2}$.\\
Similar computations to the case
$M^{2,0}\hookrightarrow\mathbb{M}^{2,1}$ give the result. Two (ISK)-spinors are needed.\\
\noindent
$\bullet$ \textit{Timelike immersion}: $\varepsilon=i$, $M^{0,2}\hookrightarrow\mathbb{M}^{0,3}$.\\
We get $\lvert\lvert
C_{0,2}\lvert\lvert^2\varphi^{\pm}=G^2_{0,2}\varphi^{\pm},$ hence
$C_{0,2}=G^2_{0,2}=0$ as the norm of $C_{0,2}$ is negative definite.
This is a similar computation to the case
$M^{2,0}\hookrightarrow\mathbb{M}^{3,0}$. Only one (ISK)-spinor is needed.\qed\\

Let us summarize these results. In the tabular below we give the number of (RSK)-(resp. (ISK)-)spinors on the surface $M^{p,q}$ solutions
of the special Killing equation \eqref{Special_Killing}, or equivalently of the Dirac equation \eqref{Dirac_equation}, which is sufficient
for the surface to be immersed, depending on the signature $(p,q)$ and on the type $\varepsilon$ of the immersion.\\\\

\begin{table}[h]
\begin{center}

{\begin{tabular}{|cc|c|c|}
\hline &   $\varepsilon$ &1& i \\
 (p,q) &   &  &\\\hline
 (0,2) &   & 2 RSK-spinor  & 1 ISK-spinor\\
 (1,1) &   & 2 RSK-spinors  & 2 ISK-spinors \\
 (2,0) &   & 1 RSK-spinor   & 2 ISK-spinors \\\hline
\end{tabular}}
\caption{\textit{Number of spinors needed}}
\end{center}
\end{table}

%We write $\varphi_1=\varphi_1^++\varphi_1^-$, then we have $\pm G\varphi_1^{\pm}=\varepsilon C\cdot\varphi_1^{\mp}$. Thus, from this relation applied twice, we obtain $G^2\varphi_1^{\pm}=||C||^2\varphi_1^{\pm}$. Similarly, we have $\pm G\varphi_1^{\pm}= -\varepsilon C\cdot\varphi_1^{\mp}$ and $G^2\varphi_2^{\pm}=||C||^2\varphi_2^{\pm}$. Here again, we need to distinguish two cases.\\
%Thus $C$ vanishes too and the Gauss and Codazzi equation are satisfied. Then, we get the conclusion by the fundamental theorem of hypersurfaces given above.$\finpreuve$

\end{document}